\documentclass[letterpaper, 10 pt, conference]{ieeeconf}

\renewcommand{\baselinestretch}{.91}

\IEEEoverridecommandlockouts                              

\usepackage{subfig}
\usepackage[top=4cm, bottom=4cm, left=2.5cm, right=2.5cm]{geometry}
\usepackage[utf8]{inputenc}
\usepackage{amsmath}
\usepackage{url}
\usepackage{graphics}
\usepackage{tabularx}
\usepackage{amsmath}
\usepackage{amssymb}

\usepackage{arydshln} 

\usepackage{booktabs} 

\usepackage{tikz}
\usepackage{pgfplots}
\usepackage{textcomp}
\usetikzlibrary{calc, shapes, arrows, decorations, positioning}
\pgfplotsset{every axis/.append style={line width=1.1pt}}
\tikzset{%
	block/.style    = {draw, rectangle, thick, minimum height = 3em,
		minimum width = {width("Spacecraft")+10pt}},
	sum/.style      = {draw, circle, node distance = 2cm}, 
	input/.style    = {coordinate}, 
	output/.style   = {coordinate} 
}

\usepackage{algorithm}
\usepackage{algorithmic}
\usepackage{amsfonts}
\usepackage{parskip}
\usepackage[outdir=./]{epstopdf}

\newcommand{\norm}[1]{\| #1 \|}

\newcommand{\lb}[1]{ {\underline {#1}}} 
\newcommand{\ub}[1]{ {\overline {#1}}} 
\newcommand{\stsum}[1]{\displaystyle\sum\limits #1} 

\newcommand{\eqdef}{\triangleq}
\newcommand{\Rr}{\mathbb R}

\newcommand{\argmin}{\mathop{\rm arg\ min}\nolimits}
\newcommand{\beq}{\begin{equation}}
	\newcommand{\eeq}{\end{equation}}
\newcommand{\ba}[1]{\begin{array}{#1}}
	\newcommand{\ea}{\end{array}}
\newcommand{\bali}{\begin{align}}
	\newcommand{\eali}{\end{align}}
\newcommand{\beqar}{\begin{eqnarray}}
	\newcommand{\eeqar}{\end{eqnarray}}
\newcommand{\beqarno}{\begin{eqnarray*}}
	\newcommand{\eeqarno}{\end{eqnarray*}}
\newcommand{\bseq}{\begin{subequations}}
	\newcommand{\eseq}{\end{subequations}}
\newcommand{\bsubeq}{\begin{subequations}}
	\newcommand{\esubeq}{\end{subequations}}
\newcommand{\bsm}{\begin{smallmatrix}}
	\newcommand{\esm}{\end{smallmatrix}}
\newcommand{\bbm}{\begin{bmatrix}}
	\newcommand{\ebm}{\end{bmatrix}}

\newcommand{\f}{\phi}
\renewcommand{\t}{\theta}
\newcommand{\p}{\psi}
\newcommand{\w}{\omega}
\renewcommand{\a}{\alpha}

\DeclareMathOperator*{\minimize}{minimize}

\newcommand{\stt}{\rm subject\ to}

\usepackage{float} 
\newfloat{algorithmfloat}{h}{loa}[section]
\floatname{algorithmfloat}{Algorithm}

\newcounter{lineofcodei}%
\newcounter{lineofcodeii}[lineofcodei]%
\newcounter{lineofcodeiii}[lineofcodeii]%
\newcounter{lineofcodeiv}[lineofcodeiii]%
\newcounter{lineofcodev}[lineofcodeiv]%
\newcounter{lineofcodevi}[lineofcodev]%
\renewcommand{\thelineofcodei}{\arabic{lineofcodei}.}%
\renewcommand{\thelineofcodeii}{\thelineofcodei\arabic{lineofcodeii}.}%
\renewcommand{\thelineofcodeiii}{\thelineofcodeii\arabic{lineofcodeiii}.}%
\renewcommand{\thelineofcodeiv}{\thelineofcodeiii\arabic{lineofcodeiv}.}%
\renewcommand{\thelineofcodev}{\thelineofcodeiv\arabic{lineofcodev}.}%
\renewcommand{\thelineofcodevi}{\thelineofcodev\arabic{lineofcodevi}.}%
\newcommand{\lineofcodeformat}{
	\settowidth{\labelwidth}{1.1.1.}%
	\setlength{\leftmargin}{\labelwidth}%
	\addtolength{\leftmargin}{2mm}%
	\setlength{\parsep}{0mm}\setlength{\itemsep}{0mm}\setlength{\rightmargin}{0mm}%
	\setlength{\topsep}{0mm}\setlength{\parskip}{0mm}\setlength{\partopsep}{0mm}}%
\newcommand{\bindi}{\begin{list}{\thelineofcodei}{\usecounter{lineofcodei}\lineofcodeformat}%
		\renewcommand{\bind}{\bindii}}%
	\newcommand{\bindii}{\begin{list}{\thelineofcodeii}{\usecounter{lineofcodeii}\lineofcodeformat}%
			\renewcommand{\bind}{\bindiii}}%
		\newcommand{\bindiii}{\begin{list}{\thelineofcodeiii}{\usecounter{lineofcodeiii}\lineofcodeformat}%
				\renewcommand{\bind}{\bindiv}}%
			\newcommand{\bindiv}{\begin{list}{\thelineofcodeiv}{\usecounter{lineofcodeiv}\lineofcodeformat}%
					\renewcommand{\bind}{\bindv}}%
				\newcommand{\bindv}{\begin{list}{\thelineofcodev}{\usecounter{lineofcodev}\lineofcodeformat}%
						\renewcommand{\bind}{\bindvi}}%
					\newcommand{\bindvi}{\begin{list}{\thelineofcodevi}{\usecounter{lineofcodevi}\lineofcodeformat}%
							\renewcommand{\bind}{\bindvii}}%
						\newcommand{\bind}{\bindi}%
						\newcommand{\eind}{\end{list}}%
\title{\LARGE \bf
Fixed-Point Constrained Model Predictive Control \\of Spacecraft Attitude
}

\author{Alberto Guiggiani, Ilya Kolmanovsky, Panagiotis Patrinos, Alberto Bemporad%
\thanks{A. Guiggiani, P. Patrinos and A. Bemporad are with IMT Institute for Advanced Studies, Lucca, Italy. At the time of submission, A. Guiggiani was a visiting scholar at the Department of Aerospace Engineering, University of Michigan, Ann Arbor, MI 48109 \texttt{\{alberto.guiggiani, panagiotis.patrinos, alberto.bemporad\}@imtlucca.it}. 
I. Kolmanovsky is with the Department of Aerospace Engineering, University of Michigan, Ann Arbor, MI 48109 \texttt{ilya@umich.edu}. }%
}

\begin{document}

\maketitle
\thispagestyle{empty}
\pagestyle{empty}

\begin{abstract}
The paper develops a Model Predictive Controller for constrained control of spacecraft attitude with reaction wheel actuators.
The controller exploits a special formulation of the cost with the reference governor like term, a low complexity addition of integral action to guarantee offset-free tracking of attitude set points, and an online optimization algorithm for the solution of the Quadratic Programming problem which is tailored to run in fixed-point arithmetic. Simulations on a nonlinear spacecraft model demonstrate that the MPC controller achieves good tracking performance while satisfying reaction wheel torque constraints.  The controller also has low computational complexity and is suitable for implementation in spacecrafts with fixed-point processors.
\end{abstract}

\section{Introduction} \label{sec:intro}

This paper addresses the development of a Model Predictive Controller (MPC)
to perform spacecraft constrained reorientation maneuvers.
The spacecraft is assumed to be actuated with reaction wheel array torques.  MPC represents an attractive framework for spacecraft attitude control given that it can deal effectively with the limited actuation authority of the reaction wheels.

Prior publications on MPC control of spacecraft attitude include \cite{Hegrenaes2005}, \cite{Silani2005} and \cite{Kalabic2014}.
In \cite{Hegrenaes2005} an explicit MPC solution is derived based on a linearized spacecraft model.  In \cite{Silani2005}, MPC is exploited for spacecraft attitude control using magnetic actuators.
In \cite{Kalabic2014}, a nonlinear MPC approach on SO(3) that uses the Lie group variational integrator-based discrete-time model is developed and shown to achieve global stabilization with respect to the spacecraft orientation.

 In the present paper we develop an online optimization-based MPC  suitable for implementation in a  processor that relies on fixed-point arithmetic. This solution is attractive for spacecraft attitude control for several reasons. Firstly, fixed-point processors are common in spacecrafts, especially in small satellites.  Secondly, while explicit MPC solutions may be implemented in fixed-point easier, the online MPC can handle changes in the model or in the constraints, as well as failure modes of the reaction wheels. Consequently approaches to online implementation of MPC are of interest. Finally, as compared to the explicit MPC solutions, online MPC implementation also requires smaller ROM size,  which is an advantage since ROM may be  limited on-board of the spacecraft, e.g., due to the need to be radiation-hardened. The challenge of embedding MPC on fixed-point processors has stimulated many recent publications (see, e.g., \cite{Longo:2013uj}, \cite{Jerez:2013va}, \cite{Guiggiani2014}); however, issues related to MPC implementation for spacecraft attitude control in a fixed-point processor have not been previously addressed.

In addition to demonstrating the fixed-point MPC solution, we also present a  formulation of the MPC optimization problem with a reference governor like term, and show that it leads to an increased constrained domain of attraction for the closed-loop system. Furthermore, a novel approach  to ensure offset-free tracking for constant attitude commands by augmenting the MPC controller with an integral action on the references will be demonstrated.  This approach is advantageous as it does not require an increase in the  complexity of MPC controller.  Through the simulations on the nonlinear spacecraft model, we  demonstrate that the resulting fixed-point MPC controller provides excellent constrained reorientation performance.

This paper is organized as follows. The spacecraft nonlinear model is described in Section \ref{sec:nonlinearModel}, followed by the problem formulation in Section \ref{sec:Motivation}.  Then, the linearized control model is defined (Section \ref{sec:controlModel}). The MPC problem setup, with a reference governor like term and integral action augmented through references, is presented in Section \ref{sec:modifiedMPC}.  Finally, the fixed-point QP solver is discussed (Section \ref{sec:fixedPoint}), its computational complexity is quantified (Section \ref{sec:ComputationalComplexity}), and closed-loop simulation results are reported (Section \ref{sec:simulations}).

\section{Spacecraft Nonlinear Model} \label{sec:nonlinearModel}
The rotational kinematics and dynamics equations of a spacecraft, considering a principal body frame fixed to its center of mass, are given by
\begin{equation}\label{eq:nlKinematics}
\begin{bmatrix} \dot{\f}(t) \\ \dot{\t}(t) \\ \dot{\p}(t)\end{bmatrix} = \tfrac{1}{c(\t)} \begin{bmatrix} c(\t) & s(\f)s(\t) & c(\f)s(\t) \\ 0 & c(\f)c(\t) & -s(\f)c(\t) \\ 0 & s(\f) & c(\f) \end{bmatrix}\begin{bmatrix}\w_1(t) \\ \w_2(t) \\ \w_3(t) \end{bmatrix}
\end{equation}
and
\begin{align} \label{eq:nlDynamics}
J_1\dot{\w}_1 &= (J_2 - J_3)\w_2\w_3+ M_1 \simeq M_1 , \nonumber \\
J_2\dot{\w}_2 &= (J_3 - J_1)\w_1\w_3 + M_2 \simeq M_2, \nonumber \\
J_3\dot{\w}_3 &= (J_1 - J_2)\w_1\w_2 + M_3 \simeq M_3,
\end{align}
where $c(\cdot) \eqdef \cos(\cdot)$ and $s(\cdot) \eqdef \sin(\cdot)$;  $\f(t)$, $\t(t)$, $\p(t)$ ($rad$) are the spacecraft roll, pitch and yaw angles, respectively; for $i = (1,2,3)$, $\omega_i(t)$ ($rad/s$) are the angular velocities, $J_i$ ($kgm^2$) are the principal moments of inertia, and $M_i$ ($Nm$) are the spacecraft moments.

We suppose that the spacecraft is equipped with 3 reaction wheels along each of its body frame axes. We consider those wheels as perfect discs with moments of inertia $\tilde{J}_i,\,i=1,2,3$, generating torques about the respective principal axes. The equations linking the spacecraft moments to the reaction wheels dynamics are defined as follows:
\begin{align} \label{eq:wheelDynamics}
M_1 &= -\tilde{J}_1 \left( \dot{\w}_1 + \ddot{\a}_1 + \dot{\a}_3\w_2 - \dot{\a}_2\w_3 \right) \simeq -\tilde{J}_1(\dot{\w}_1 + \ddot{\a}_1), \nonumber \\
M_2 &= -\tilde{J}_2 \left( \dot{\w}_2 + \ddot{\a}_2 + \dot{\a}_1\w_3 - \dot{\a}_3\w_1 \right)\simeq -\tilde{J}_2(\dot{\w}_2 + \ddot{\a}_2), \nonumber \\
M_3 &= -\tilde{J}_3 \left( \dot{\w}_3 + \ddot{\a}_3 + \dot{\a}_2\w_1 - \dot{\a}_1\w_2 \right)\simeq -\tilde{J}_3(\dot{\w}_3 + \ddot{\a}_3),
\end{align}
where $\dot{\a}_i \,(rad/s)$ are the wheels rotational speeds, and $\ddot{\alpha}_i (rad/s^2)$ are the wheels accelerations.

Putting together equations \eqref{eq:nlKinematics}-\eqref{eq:wheelDynamics} we can formulate a ninth-order state-space nonlinear model with state $x = \begin{bmatrix} \f & \t & \p & \w_i & \dot{\a}_i \end{bmatrix}'$, for $i=1,2,3$, that evolves according the following ODEs (dependencies on time are omitted):

\begin{subequations} \label{eq:totalSystem}
\begin{align} 
\dot{\f} &= \tfrac{1}{c(\t)}\left( c(\t)\w_1 + s(\f)s(\t)\w_2 + c(\f)s(\t)\w_3 \right)  \\
\dot{\t} &= \tfrac{1}{c(\t)}\left( c(\f)c(\t)\w_2 - s(\f)s(\t)\w_3 \right)  \\
\dot{\p} &= \tfrac{1}{c(\t)}\left( s(\f)\w_2 + c(\f)\w_3 \right)  \\
\dot{\w}_1 &= \tfrac{1}{J_1 + \tilde{J}_1}\left( (J_2-J_3)\w_2\w_3 - \tilde{J}_1(\dot{\a}_3\w_2 - \dot{\a}_2\w_3) - u_1 \right)  \\
\dot{\w}_2 &= \tfrac{1}{J_2 + \tilde{J}_2}\left( (J_3-J_1)\w_1\w_3 - \tilde{J}_2(\dot{\a}_1\w_3 - \dot{\a}_3\w_1) - u_2 \right)  \\
\dot{\w}_3 &= \tfrac{1}{J_3 + \tilde{J}_3}\left( (J_1-J_2)\w_1\w_2 - \tilde{J}_3(\dot{\a}_2\w_1 - \dot{\a}_1\w_2) - u_3 \right)  \\
\ddot{\a}_i &= \tfrac{1}{\tilde{J}_i}u_i,\; i=1,2,3,
\end{align}
\end{subequations}
where $u_1$, $u_2$, $u_3$ are the torques exerted on the wheels by the electric actuators.

\section{Motivation and Contribution}\label{sec:Motivation}
The control objective is to track a reference spacecraft orientation
\begin{equation} \label{eq:regulationProblem}
\begin{bmatrix} \f(t) \\ \t(t) \\ \p(t)  \end{bmatrix} \rightarrow \begin{bmatrix} \f^r(t) \\ \t^r(t) \\ \p^r(t)  \end{bmatrix} = r(t)
\end{equation}
subject to polytopic constraints in the form
\begin{equation} \label{eq:constraints}
\lb{z} \leq z_c \leq \ub{z}, \; z_c = C_c\begin{bmatrix}  \f(t) \\ \t(t) \\ \p(t) \\ \dot{\a}_i \end{bmatrix}+D_cu(t),
\end{equation}
for $i=1,2,3$. In other words, the control framework must be able to handle constraints on the spacecraft orientation, the reaction wheel speeds, and the control inputs.

To this end, we formulate a Model Predictive Control setup that: (a) is able to solve problem \eqref{eq:regulationProblem}-\eqref{eq:constraints} taking into account both spacecraft and reaction wheels dynamics; (b) minimizes the computational impact and memory requirements; (c) can be deployed on fixed-point microcontrollers or as an embedded MPC-on-a-chip device. This is achieved in the following way. 

First, we define a reduced-order control model, halving the overall problem size while retaining a description of the significant spacecraft and wheels dynamics (Section \ref{sec:controlModel}). We propose a modified MPC formulation with virtual optimization variables and references that guarantees offset-free tracking while lowering the prediction horizon requirements, further reducing the QP problem size (Section \ref{sec:modifiedMPC}). Finally, we assign the solution of the QP problem to an algorithm tailored for execution on embedded devices, exploiting fixed-point arithmetics to reduce computational load and memory footprint (Section \ref{sec:fixedPoint}). 

\section{Control model} \label{sec:controlModel}

The nonlinear model \eqref{eq:totalSystem} is too complex as a prediction model for an embedded MPC controller. We need to formulate a linear, reduced-order model, that however is still capable of capturing the significant dynamics of the system. 



Since the angular momentum is conserved and neglecting all but linear terms in \eqref{eq:nlDynamics}-\eqref{eq:wheelDynamics}, we set the linear model
\begin{equation}\label{eq:relationWheelsSat}
\dot{\w}_i = -\tfrac{\tilde{J}_i}{J_i}\ddot{\a}_i,\; i=1,2,3,
\end{equation}
hence
\begin{equation}
\ddot{\a}_i = \tfrac{J_i}{J_i\tilde{J}_i-\tilde{J}_i^2}u_i \eqdef \mathbf{f}\left(J_i,\tilde{J}_i\right)u_i, \; 1= 1,2,3.  
\end{equation}

We are now ready to define the reduced-order model for MPC design in LTI state-space form, linearized for small angles, as follows

\begin{align}\label{eq:linearizedModel}
\overbrace{\begin{bmatrix} \dot{\f} \\ \dot{\t} \\ \dot{\p} \\ \ddot{\a}_i \\ \ddot{\a}_2 \\ \ddot{\a}_3  \end{bmatrix} }^{\dot{x}} &= A_C \overbrace{\begin{bmatrix} \f \\ \t \\ \p \\ \dot{\a}_1 \\ \dot{a}_2 \\ \dot{a}_3  \end{bmatrix} }^{x} + B_C \overbrace{ \begin{bmatrix} u_1 \\ u_2 \\ u_3 \end{bmatrix} }^{u}, \\
A_C &\eqdef \begin{bmatrix} \begin{array}{c|c} \mathbf{0}^{3\times 3} & \begin{matrix} -\tfrac{\tilde{J}_1}{J_1} & 0 & 0 \\ 0 & -\tfrac{\tilde{J}_2}{J_2} & 0 \\ 0 & 0 & -\tfrac{\tilde{J}_3}{J_3} \end{matrix} \\ \hline \mathbf{0}^{3\times 3} & \mathbf{0}^{3\times 3} \end{array}\end{bmatrix}, \nonumber \\
B_C &\eqdef \begin{bmatrix} \begin{array}{c} \mathbf{0}^{3\times 3} \\ \hline \begin{matrix} \mathbf{f}\left(J_1,\tilde{J}_1\right) & 0 & 0 \\ 0 & \mathbf{f}\left(J_2,\tilde{J}_2\right) & 0 \\ 0 & 0 & \mathbf{f}\left(J_3,\tilde{J}_3\right) \end{matrix}  \end{array}\end{bmatrix}. \nonumber
\end{align}


\section{MPC formulation} \label{sec:modifiedMPC}

While designing an MPC scheme, the choice of an appropriate prediction horizon is a critical step to ensure proper controller performance. The prediction horizon should ideally cover the time needed to perform the required maneuvers (e.g., on a reference step change, the controller should be able to ``see" in the future when the orientation has approached the new reference). Therefore, one should know \emph{a priori} bounds on the reference variations in order to select a prediction horizon large enough to account for them. 

Moreover, we highlight two more issues: (a) the spacecraft maneuvers are usually ``slow" compared with the controller sampling rates, leading to requirements for large prediction horizon and consequently a large size of the Quadratic Program (QP) associated with the MPC controller; and (b) a prediction horizon tailored to account for the upper bound on the reference variations may cause the smaller maneuvers to be excessively slow (in case longer horizons lead to less aggressive control actions).

To address these issues, we propose a modified MPC setup with additional virtual optimization variables and references. The approach is similar to the one in \cite{limon_mpc_2008}, and can be described as follows.

We define a new optimization vector, where the variations on the control input (a standard choice for MPC tracking problems) are extended with a new vector of length equal to the number of system states, becoming $\begin{bmatrix} \Delta u_k' & \tilde{x}_k' \end{bmatrix}'$, where $\Delta u_k \in \Rr^3$ and $\tilde{x}_k \in \Rr^6$, $k = 0, ..., N-1$, and set the new cost function
\begin{align}\label{eq:newCost}
V\left(x, \tilde{r}, \Delta u, \tilde{x} \right) &= \norm{x_N - \tilde{r}}_{P_f} + \\
\stsum_{k=0}^{N-1}&{ \norm{(x_k - \tilde{x}_k) - \tilde{r}}_{Q_1} + \norm{\tilde{x}_k}_{Q_2} + \norm{\Delta u_k}_R }, \nonumber
\end{align}
where $N$ is the prediction horizon, $\tilde{r}$ is the reference signal for spacecraft attitude and wheels speeds, $\tilde{x} = \begin{bmatrix} \tilde{x}_0' & \hdots & \tilde{x}_{N-1}' \end{bmatrix}'$, $\Delta{u} = \begin{bmatrix} \Delta{u}_0' & \hdots & \Delta{u}_{N-1}' \end{bmatrix}'$.

The optimal control problem to be solved at each sampling step becomes

\begin{align} \label{eq:modifiedProblem}
V^*\left(x(t), \tilde{r}(t), u(t)\right) &= \min_{\Delta u, \tilde{x}} V\left(x(t), \tilde{r}(t), \Delta u, \tilde{x}\right) \\
\stt\quad  x_{k+1} &= Ax_k + Bu_k, \; k = 0, ..., N-1 \nonumber \\
u_k &= u_{k-1} + \Delta u_k, \; k = 0, ..., N_c-1 \nonumber \\
u_k &= u_{k-1}, \; k = N_c, ..., N-1 \nonumber \\
u_{-1} &= u(t), \; x_0 = x(t) \nonumber \\
(x_k, u_k) &\in \mathcal{Z}, \; k = 0, ..., N-1, \nonumber
\end{align}

where $N_c$ is the control horizon and $\mathcal{Z}$ is the polytope defined in \eqref{eq:constraints}.

Looking at the cost function \eqref{eq:newCost}, notice that there is no direct penalty for the state being far from the reference, as in standard tracking MPC. Instead, the optimizer has the freedom to trade-off between the discrepancy state/reference and the magnitude of the auxiliary variables $\tilde{x}$, for each prediction step. This trade-off is also influenced by properly tuning the weight matrices $Q_1$, $Q_2$, $P_f$.

The practical impact on the control performance when adopting the modified MPC setup \eqref{eq:newCost}-\eqref{eq:modifiedProblem} is that the prediction horizon can be chosen independently of the duration of the spacecraft maneuvers, and that the simulated closed-loop system shows stable behavior for significantly smaller prediction horizons. Similar effects have been observed in the orbital maneuvering study \cite{weiss_model_2012}.

\begin{figure}
\centerline{\includegraphics[width=\columnwidth]{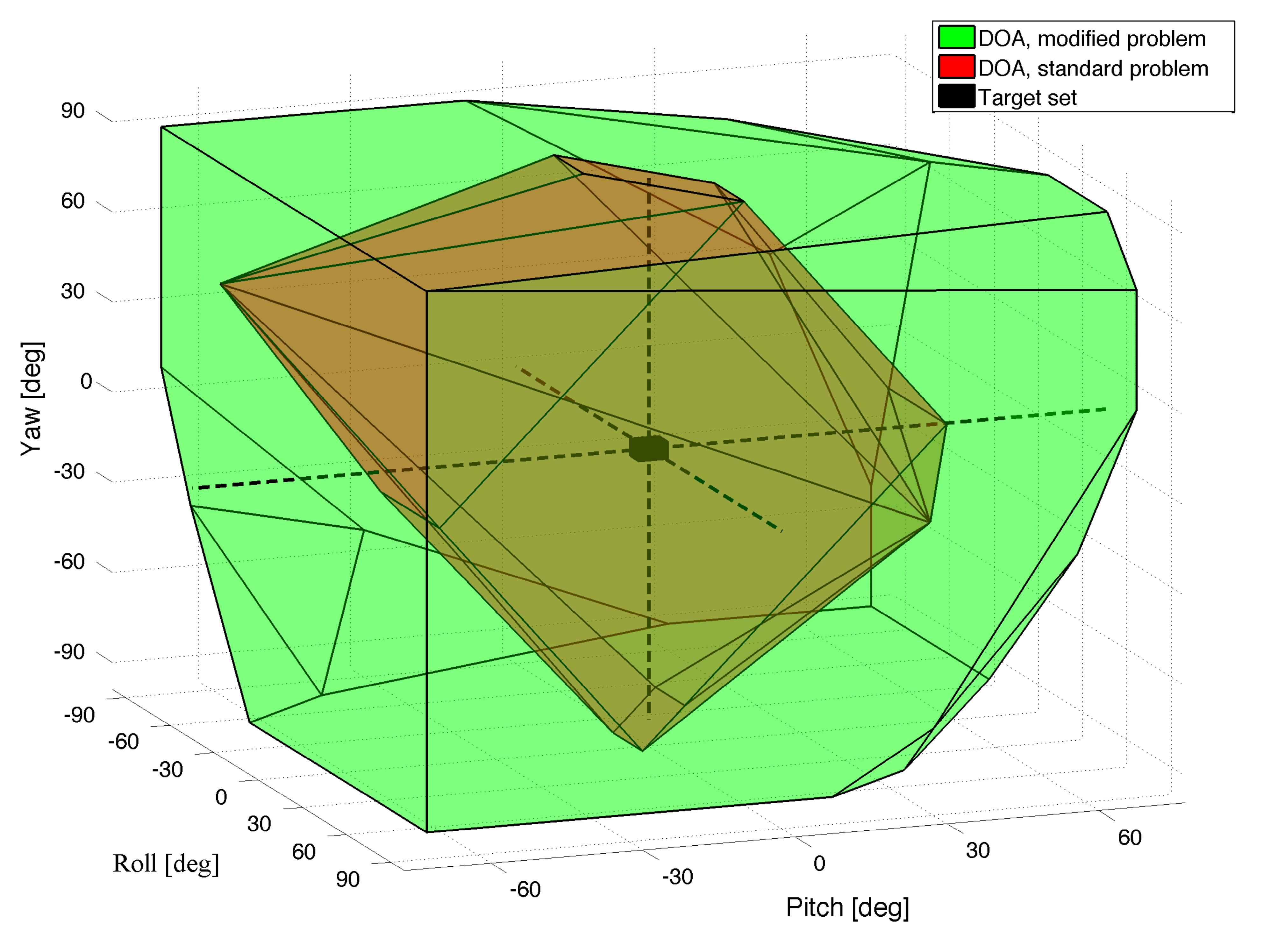}}
\vspace{-2mm}
\caption{Domain of attraction for the modified problem \eqref{eq:newCost}-\eqref{eq:modifiedProblem} (green) compared to the domain of attraction for a standard MPC setup (brown).}
\label{fig:doas}
\end{figure}

Figure \ref{fig:doas} shows the domain of initial conditions for which the controller is able to drive the state trajectories to a target set (black) of radius $0.05 \;rad$, for a spacecraft simulated with the full nonlinear model \eqref{eq:totalSystem}. The domains of attraction are computed by performing a grid search on the state space and running closed loop simulations with two controllers. The green set is obtained when the controller solves the modified problem proposed in \eqref{eq:modifiedProblem}, the red set when is instead solved a standard MPC tracking problem (i.e., removing the auxiliary optimization variables $\tilde{x}$), while leaving all the other parameters unchanged. Sampling time is $0.5s$, and the prediction and control horizons are set to $10$ and $2$, respectively. Wheel speeds are constrained in $[-10, 10] \, rad/s$, and control moments in $[-1, +1] \, Nm$. Spacecraft moments of inertia are $J_1=3000\,kgm^2$, $J_2=1500\,kgm^2$, and $J_3=2000\,kgm^2$, while wheels moments of inertia are $\tilde{J}_i = 50 \, kgm^2,\; i=1,2,3$.
The weights in the cost function are 
$$Q_1 = \begin{bmatrix} \begin{array}{c|c} 100\cdot \mathbf{I}^{3} & \mathbf{0}^{3\times 3} \\ \hline \mathbf{0}^{3\times 3} & 0.1\cdot \mathbf{I}^{3} \end{array}\end{bmatrix},$$ $$Q_2 = 50\cdot \mathbf{I}^{6}, R = 0.01\cdot \mathbf{I}^{3},$$
where $\mathbf{I}^{n}$ is the identity matrix of size $n$. The terminal weight $P_f$ is the solution of the Riccati equation associated to the LQR problem.

Simulation results show how, even with a prediction horizon $N = 10$ (that gives the controller a prediction window of $5$ seconds only), we obtain a domain of attraction spanning up to $[-\pi/2, \pi/2]$ radians for the roll and yaw angles, and up to $[-\pi/3, \pi/3]$ radians for the pitch angle. This is a remarkable result, considering that the controller operates with a prediction model linearized for small angles.

However, the control development is not complete yet, as it does not guarantee offset-free tracking in presence of model uncertainties, and is not able to reject any constant external disturbance, such as caused by gravity gradients or solar radiation pressure. To address this problem, we add an integral action on the tracking error.

There are several ways to achieve integral action in a controller. For example, one can extend the system state with an integrator state, or with a disturbance observer. However, we aim at keeping the QP size as small as possible for efficient embedded computations, and therefore avoid increasing the model size.

The proposed approach works as follows. We exploit the additional degree of freedom in the optimizer given by the virtual reference, and add the integral action to the reference itself. We separate the \emph{true} reference $\bar{r}(t)$ (i.e., the actual desired spacecraft orientation) with the \emph{actual} reference $r(t)$ fed to the controller, by means of the following equation
\begin{equation} \label{eq:integral}
r(t) = \bar{r}(t) - \int_0^t\left( \begin{bmatrix} \f(\tau) \\ \t(\tau) \\ \p(\tau) \end{bmatrix} - \bar{r}(\tau)\right) \,\mathrm{d}\tau.
\end{equation}

\begin{figure}
\centerline{\includegraphics[width=\columnwidth]{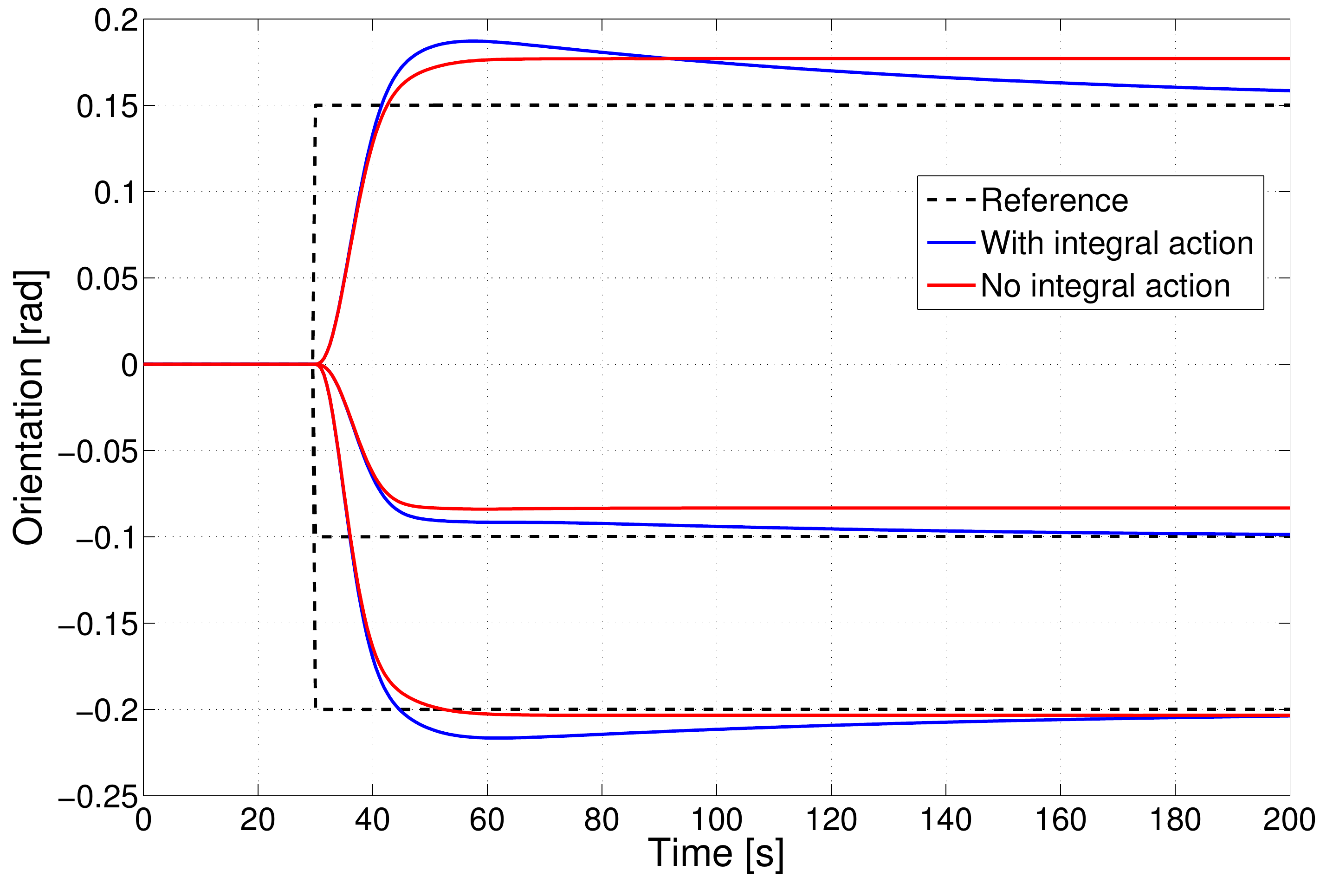}}
\vspace{-2mm}
\caption{Tracking of reference roll, pitch and yaw angles (dashed lines) using a control scheme with the proposed integral action on the reference (blue lines) and without (red lines).}
\label{fig:offsetFree}
\end{figure}

With this approach we achieve offset-free tracking. The computation of the integral action \eqref{eq:integral} can be easily performed by an external integrator without affecting the complexity of the controller. The effect is depicted in Figure \ref{fig:offsetFree}; we observe that the controller is able to reject constant disturbances thanks to the external integral action.

\section{Fixed-point QP solver} \label{sec:fixedPoint}

\subsection{Fixed-Point Computations} \label{sub:FPComputations}
By \emph{fixed-point} we refer to way to represent numbers in digital processors, meaning a mapping from real-world numbers to sequences of binary values. A fixed-point data type is denoted by four parameters: (1) the total word length $w$, (2) the signedness $s$, (3) the number $r$ of bits for the integer part, and (4) the number $p$ of bits for the fractional part.

The choice of fixed-point computations can have a great positive impact when trying to minimize computational effort, power consumption, and chip size \cite{Kerrigan}. Moreover, fixed-point representations are of particular interest in embedded applications as they grant hardware support for additions and multiplications on nearly all computing devices.

However, when a number is coded as a binary word with fixed-point representation, a fixed number $p$ of bits is assigned to its fractional part: this leads to a limited precision equal to $2^{-(p+1)}$, defined as the difference between two successive values, and therefore to round-off errors. In details, multiplying two fixed-point numbers $\zeta=\gamma\xi$ leads to the fixed-point representation $\texttt{fi}(\zeta)$ of $\zeta$, with $|\zeta-\texttt{fi}(\zeta)|\leq 2^{-(p+1)}$. In case of inner product operation between two vectors $x,y\in\Rr^n$, the total error can be bounded as
\begin{equation}
|x'y-\texttt{fi}(x'y)|\leq 2^{-(p+1)}n.
\end{equation}
Moreover, the fixed number $r$ of bits for the integer part leads to the occurrence of overflow errors when trying to represent values outside the admissible range $[-2^r, 2^r-1]$.

When solving problem \eqref{eq:modifiedProblem} with fixed-point arithmetics, it is mandatory that the underlying algorithm: (a) guarantees convergence to an acceptable solution despite round-off errors; and (b) prevents the occurrence of overflow errors. An algorithm that satisfies these conditions was introduced in \cite{Patrinos:2013to} as a modified version of the \emph{GPAD} algorithm (cf. \cite{GPAD_TAC}), and will be briefly summarized in the following section.

\subsection{Fixed-Point Dual Gradient Projection Algorithm ({FP-GPD})} \label{sub:FP-GPD}

By means of simple algebraic steps it is possible to formulate \eqref{eq:modifiedProblem} in a condensed Quadratic Programming (QP) form (cf. \cite{Maciejowski:2002wd}), obtaining
\begin{align}\label{eq:optimizationProblemCondensed}
\minimize\quad& V(z) = \tfrac{1}{2}z'Hz+ h'z\\
\stt\quad& g(z)\leq 0, \nonumber
\end{align}
where $V:\Rr^n\to\Rr$ is differentiable and strongly convex, $z\in\Rr^n$ denotes the vector of optimization variables, $H\in \Rr^{n\times n}$ is a positive definite matrix, and $g(z)=Dz-d$ is an affine mapping.

Consider an approximate solution of problem \eqref{eq:optimizationProblemCondensed} as follows. A vector $z\in\mathbb{R}^{n}$ is an $\left(\varepsilon_{V},\varepsilon_{g}\right)$-optimal solution for two nonnegative constants $\varepsilon_{V}$ and $\varepsilon_{g}$ if 
\begin{align}\label{epssolution}
 V(z)-V^{*}&\leq\varepsilon_{V},\\
 \left\Vert \left[g(z)\right]_{+}\right\Vert_\infty&\leq\varepsilon_{g},\label{epssolution2}
 \end{align}\\ 
 where $\left[\cdot\right]_{+}$ denotes the Euclidean projection on the nonnegative orthant and $\left\Vert\cdot\right\Vert_\infty$ is the infinity norm.\\
Starting from \eqref{eq:optimizationProblemCondensed} and relaxing the inequality constraints, the \emph{dual problem} is:
 \begin{align}\label{eq:dualProblem}
\mathbb{D}:&&\Phi ^{*}=\underset{y\geq0}{\max}\left(\underset{z\in\Rr^n}{\min}\;V(z)+y'g(z)\right)
\end{align}
As it is well known, for QPs strong duality holds provided that~\eqref{eq:optimizationProblemCondensed} is feasible, i.e., $V^{\star}=\Phi^{\star}$ and the unique primal solution can be recovered by any dual optimal solution via $z^\star=-H^{-1}(D'y+h)$.

The main goal is to apply the optimization algorithm to the dual problem \eqref{eq:dualProblem} and converge to a primal solution with given suboptimality and infeasibility, as in \eqref{epssolution} and \eqref{epssolution2}. 

The proposed Gradient Projection algorithm applied to problem \eqref{eq:dualProblem} performs the following iterations,
\begin{subequations}\label{eq:DualGP}
\begin{align}
z_{(\nu)}&=\underset{z\in\Rr^n}{\argmin}\left(V(z)+y_{(\nu)}'g(z)\right),\label{eq:DualGPa}\\
y_{(\nu+1)}&=\left[y_{(\nu)}+\tfrac{1}{L_{\Phi}}g(z_{(\nu)})\right]_+,\label{eq:DualGPb}
\end{align}
\end{subequations}
where $L_{\Phi}$ is the Lipschitz constant of the dual problem Hessian. The algorithm is initialized with $y_{(0)} = 0$, and stopped as soon as conditions \eqref{epssolution}-\eqref{epssolution2} are met.

In practical implementations, the most expensive computation needed for steps \eqref{eq:DualGPa} and \eqref{eq:DualGPb} is one matrix-vector product each, with matrices sizes equal to the number of primal and dual variables. The whole algorithm is division-free, and does not require matrix inversion procedures as in most second-order solvers. Those are key strengths that make the algorithm appealing for embedded implementation setups.


Moreover, the algorithm is robust to finite-precision computations, with proven convergence properties when executed in fixed-point arithmetic. Specifically, it is possible to (a) define lower bounds on the number $p$ of fractional bits such that a desired solution quality \eqref{epssolution}-\eqref{epssolution2} is achieved, and (b) define lower bounds on the number $r$ of integer bits such that avoidance of overflow error is guaranteed. We refer the reader to \cite{Patrinos:2013to} for a detailed analysis of those results.

%
%

\section{Computational complexity} \label{sec:ComputationalComplexity}

In this section we detail the computational complexity of the predictive controller based on the MPC formulation introduced in Section \ref{sec:modifiedMPC}, supported by the fixed-point QP solver proposed in Section \ref{sec:fixedPoint}.

The QP algorithm has been implemented in library-free ANSI-C for a streamlined deployment on hardware platforms. The resulting code is composed by two functions only; an \texttt{init} function, to be called once, which initializes the problem data, and a \texttt{step} function, to be called at each sampling step, which accepts spacecraft orientation measurements and references as input, and outputs the control signals for the reaction wheels torques.

Table \ref{tab:complexity} analyzes three different controller configurations: the first with input constraints only (upper and lower bounds on all the control signals); the second adding constraint on the input variations; the last including also constraints on both spacecraft orientation (inclusion zone constraints) and wheels speed. The \emph{Variables} column reports the number of primal and dual variables of the resulting QP problem. The \emph{Size} column shows the memory requirements to store the problem data and to execute the code functions. Finally, the \emph{Oper./Iter} column reports the number of fixed-point operations (multiplications and additions) required to complete an algorithm iteration \eqref{eq:DualGP}.

Results show how, thanks to the modified MPC formulation and the choice of the fixed-point QP solver, we are able to maintain the QP problem small and solve it efficiently with minimal memory footprint and computational burden. Moreover, given specific hardware it is possible to estimate precisely the time required to perform a single iteration (since the algorithm performs linear-algebra computations on matrices of pre-determined sizes). Then, given the sampling time, one can determine the number of iterations that can be computed within the sampling period. Finally, using the results in \cite{rubagotti2014stabilizing} one can formulate the MPC problem with stability and recursive feasibility guarantees.

\begin{table}
\caption{Controller Complexity}
\label{tab:complexity}
\begin{center}
\begin{tabular}{|c||c|c|c|}
\hline
\textbf{Constraints} & \textbf{Variables} & \textbf{Size} [$KB$] & \textbf{Oper./Iter}\\
\hline
\hline
$u$ & $18\, |\, 12$ & $8.7\, |\, 4.2$ & $460$ \\
$u + \Delta u$ & $18\, |\, 24$ & $8.7\, |\, 4.5$ & $900$ \\
$u + \Delta u + x$ & $\,\,\,18\, |\, 144$ & $37.5\, |\, 7.7\,\,\,$ & $5300$ \\
\hline
\end{tabular}
\end{center}
\end{table}

\section{Simulations}\label{sec:simulations}

The following simulations are aimed to investigate the closed-loop behavior when a controller based on the proposed MPC setup is connected to a spacecraft simulated with the nonlinear model of Section \ref{sec:nonlinearModel}.

\subsection{Sinusoidal References Tracking} \label{sub:sinusoidal}

The spacecraft is required to track sinusoidal references, with varying amplitudes and frequencies for the three angles. The MPC parameters are the same used in the design that produced Figure \ref{fig:doas}.

\begin{figure}
  \centering
  \subfloat[Control moments]{\label{fig:sinRef_inputs}\includegraphics[width=\columnwidth]{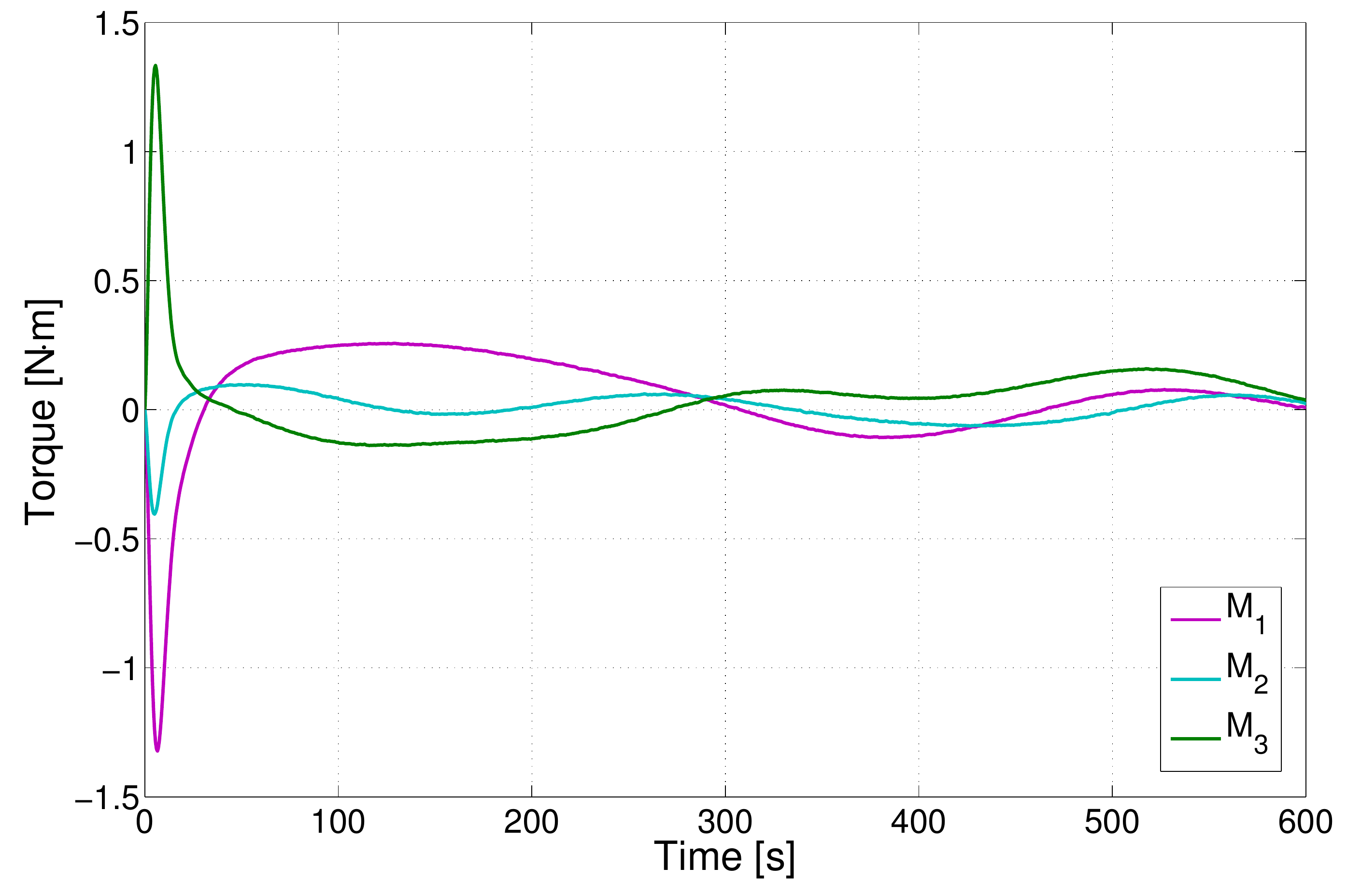}}\\              
  \subfloat[Wheels rotational speeds]{\label{fig:sinRef_wheels}\includegraphics[width=\columnwidth]{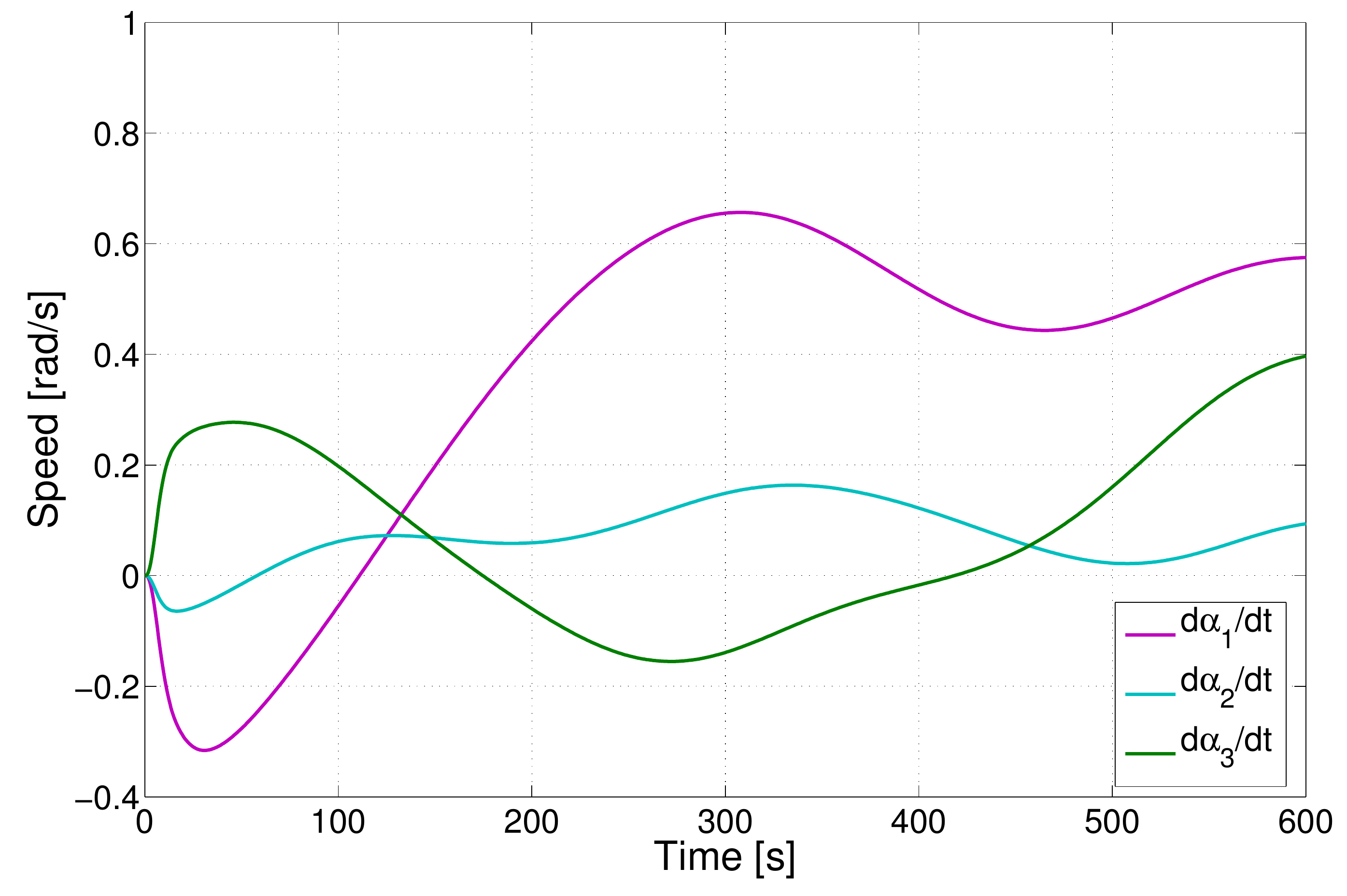}}\\
  \subfloat[Spacecraft orientation]{\label{fig:sinRef_satPos}\includegraphics[width=\columnwidth]{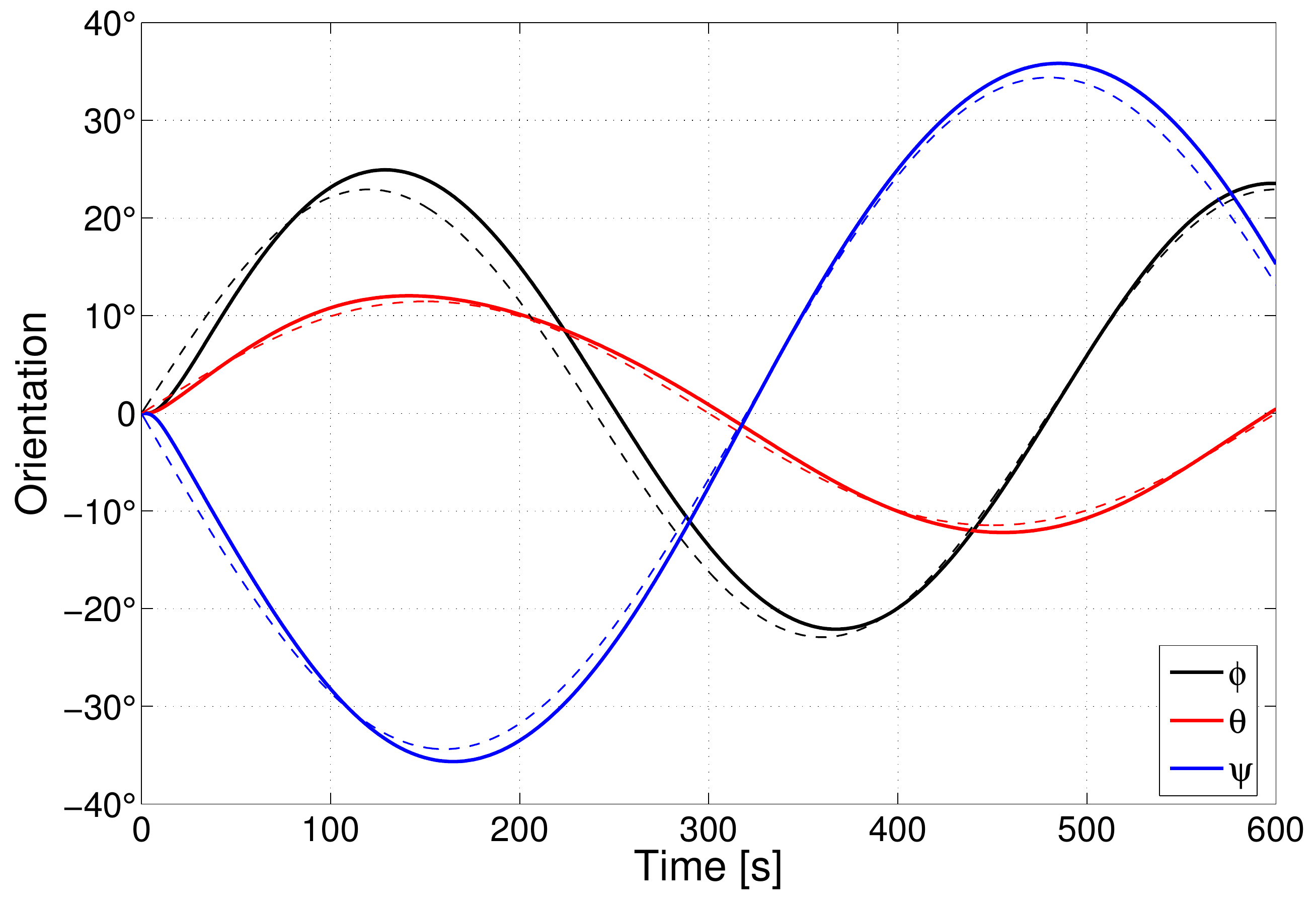}}
  \caption{Closed-loop simulation for sinusoidal reference tracking: (a) control inputs (torques applied to the three wheels); (b) reaction wheels rotational speeds; (c) spacecraft roll, pitch and yaw angles (solid lines) compared to reference trajectories (dashed lines).}
  \label{fig:sinRef}
\end{figure}

Figure \ref{fig:sinRef} shows the result of the closed-loop simulation with plots for the control inputs (\ref{fig:sinRef_inputs}), wheels rotational speeds (\ref{fig:sinRef_wheels}), and spacecraft orientation compared to reference trajectories (\ref{fig:sinRef_satPos}). 

The closed-loop behavior is consistent with the references, despite the fact that the controller relies on a reduced system model for predictions. It has to be noted that, in the current MPC formulation, the controller is not aware of future references. If such information is available, it can be incorporated into the prediction model further improving the controller performances.

\subsection{Fixed-Point Accuracy}

The goal of this simulation is to evaluate the robustness of the proposed QP solver with respect to finite-precision number representations. 

\begin{figure}
\centerline{\includegraphics[width=\columnwidth]{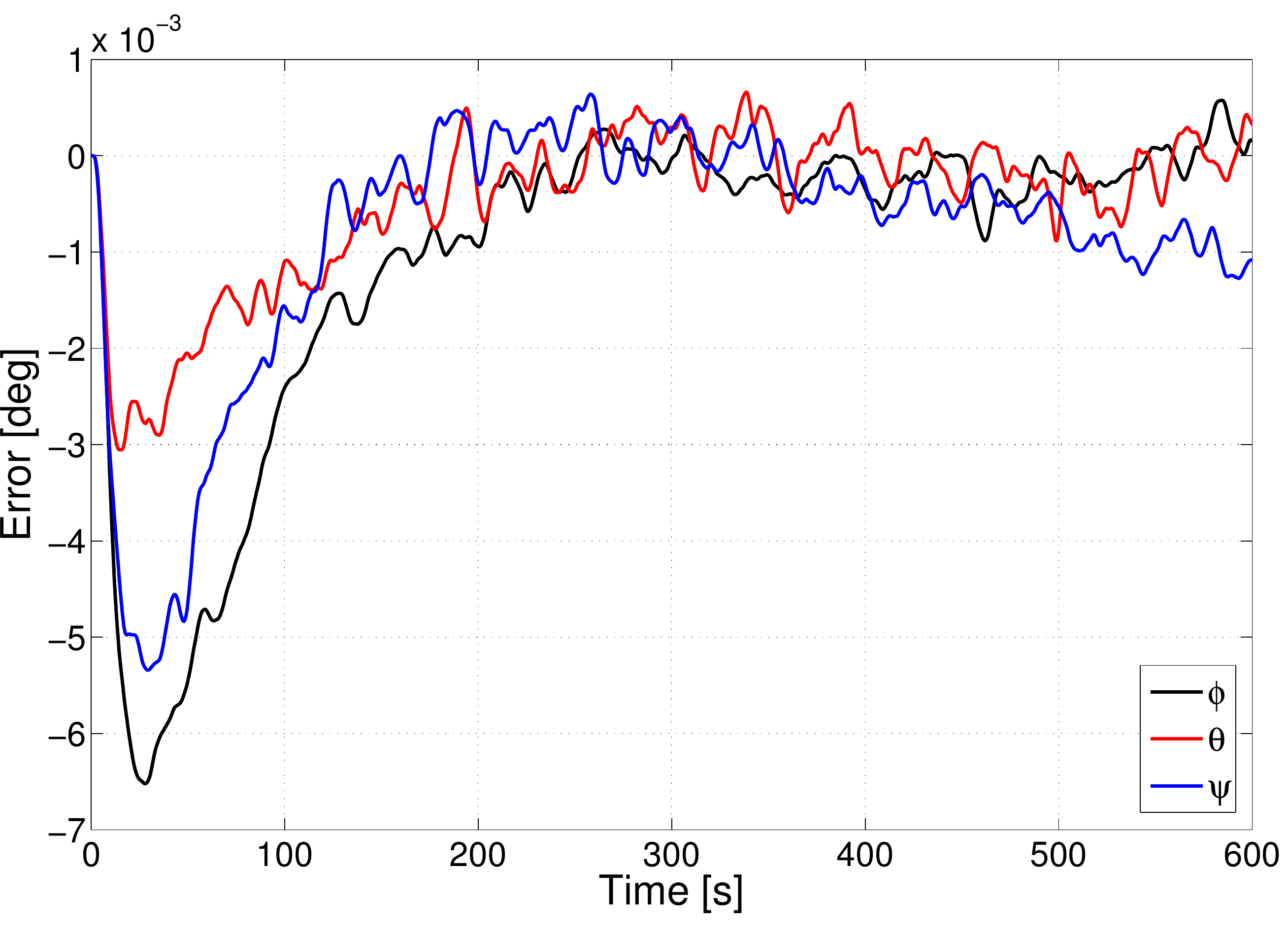}}
\vspace{-2mm}
\caption{Difference in spacecraft orientation using the fixed-point QP solver compared to 64-bit floating-point solver.}
\label{fig:error}
\end{figure}

We ran the closed-loop simulation of Section \ref{sub:sinusoidal} first with 32-bit fixed-point arithmetic, of which 16 bits for the fractional part. Then, we switched to 64-bit, double precision floating point arithmetic. The plots in Figure \ref{fig:error} show the discrepancy in the resulting spacecraft orientation. We observe how the divergence, although larger when the control action is stronger, remains in the order of $10^{-3}$ degrees for the whole maneuver.

\subsection{Rest-to-Rest Orientation Maneuver}

The purpose of the last simulation is to show the closed-loop behavior when performing rest-to-rest orientation maneuvers, for different actuator saturations.

The controller is set to track a reference step change, where $\begin{bmatrix} \f^r(t) \\ \t^r(t) \\ \p^r(t) \end{bmatrix}: \begin{bmatrix} 0 \\ 0 \\ 0 \end{bmatrix} \rightarrow \begin{bmatrix} 0.08 \\ -0.03 \\ -0.1 \end{bmatrix} [rad]$. The maneuver is repeated for loose input constraints, $u(t) \in [-3, +3] Nm$, and for tight input constraints, $u(t) \in [-0.2, +0.2] Nm$. All the other parameters are the same as in the design that generated Figure \ref{fig:doas}.

\begin{figure}
\centerline{\includegraphics[width=\columnwidth]{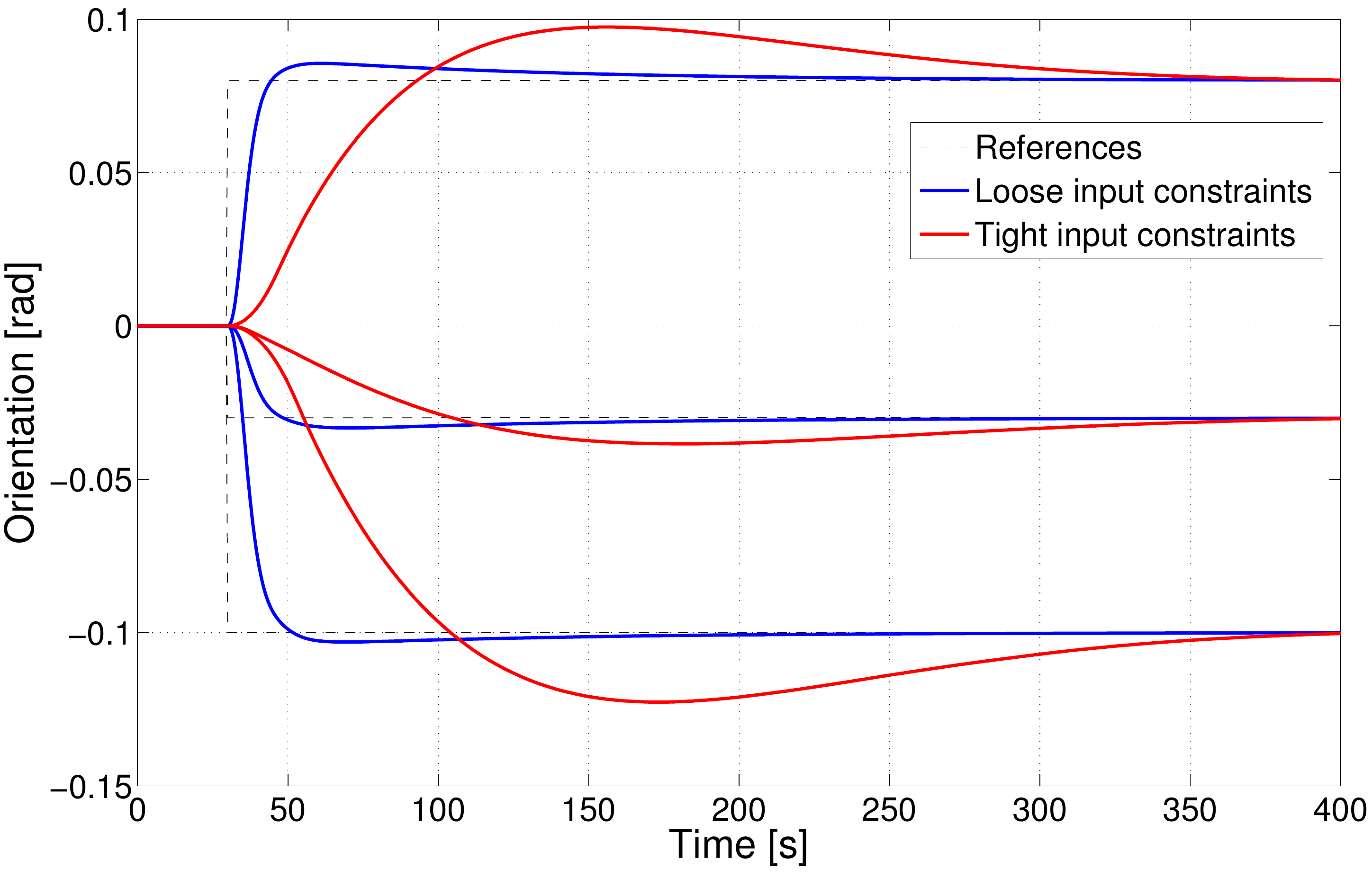}}
\vspace{-2mm}
\caption{Rest-to-rest orientation maneuver for loose (blue) and tight (red) input constraints.}
\label{fig:restToRest}
\end{figure}

Figure \ref{fig:restToRest} shows how the controller is able to react to the different input constraints and completes the rest-to-rest orientation maneuver, even in the case of strongly constrained actuators.

\section{Conclusion} \label{sec:conclusion}
This paper presented a fixed-point Model Predictive Control framework for spacecraft attitude tracking with reaction wheels actuators. A reduced control model is chosen such that the computational load is minimized while retaining the significant spacecraft and reaction wheel dynamics. Then, a modified cost function is presented, allowing for a significant increase in the controller domain of attraction and a reduction in the prediction horizon, thus keeping the QP problem size small. Moreover, an external integral action on the reference is introduced, granting offset-free tracking without adding complexity to the controller itself. Finally, the QP problem is assigned to a fixed-point solver to further reduce memory footprint and computational burden, while retaining the closed-loop performance. As a result, an efficient and lightweight ANSI-C implementation is obtained; its minimal memory and computational requirements make it suitable for deployment on low-power embedded devices for attitude control.

\bibliographystyle{IEEEtran}
\bibliography{IEEEabrv,bibliography}

\end{document}